\documentclass[a4paper]{article}
\usepackage{lmodern}
\pdfoutput=1
\usepackage{graphicx,wrapfig}
\usepackage{lipsum}
\usepackage{amsmath,amssymb,mathrsfs}
\usepackage{hyperref}
\usepackage{url}
\usepackage{mathtools}
\usepackage{changepage}
\usepackage[ruled,vlined]{algorithm2e}
\usepackage{amsmath}
\usepackage{multirow}
\usepackage{wrapfig}
\usepackage{subfig}
\usepackage{color}
\usepackage{epsfig,enumerate,amsmath,amsfonts,amssymb,amsthm,mathrsfs,ifpdf}
\usepackage{indentfirst,relsize}
\usepackage{setspace,graphicx}
\usepackage{latexsym}
\usepackage[all]{xy}
\usepackage[dvipsnames]{pstricks}
\usepackage{pst-grad} 
\usepackage{pst-plot} 
\usepackage[margin = 2.50cm]{geometry}

\usepackage{authblk} 
\usepackage{algpseudocode}
\usepackage{color,soul}
\usepackage{enumitem}

\usepackage{orcidlink}

\newcommand{\remove}[1]{}

\newtheorem{theorem}{Theorem}
\newtheorem{lemma}[theorem]{Lemma}

\newcounter{Case}[theorem]
\newtheorem{case}[Case]{Case}

\title{Vertex-edge domination on subclasses of bipartite graphs}

\author[1]{Arti Pandey\footnote{arti@iitrpr.ac.in}}
\author[2]{Kaustav Paul\footnote{kaustav.20maz0010@iitrpr.ac.in}}
\author[3]{Kamal Santra \orcidlink{0009-0006-5997-1452} \footnote{kamal.7.2013@gmail.com, kamal.santra@iitg.ac.in}\footnote{This research was carried out when the author was a Research Associate at the Indian Institute of Technology Ropar.}}
\affil[1, 2]{Department of Mathematics\\
	
	Indian Institute of Technology Ropar\\
	
	Rupnagar, 140001, Punjab, India}
	
	\affil[3]{Department of Mathematics\\
		
		Indian Institute of Technology Guwahati\\
		
		Guwahati, 781039, Assam, India}
\date{}

\begin{document}

\maketitle
\begin{abstract}
Given a simple undirected graph $G = (V, E)$, the open neighbourhood of a vertex $v \in V$ is defined as $N_G(v) = \{u \in V \mid uv \in E\}$, and the closed neighbourhood as $N_G[v] = N_G(v) \cup \{v\}$. A subset $D \subseteq V$ is called a vertex-edge dominating set if, for every edge $uv \in E$, at least one vertex from $D$ appears in $N_G[u] \cup N_G[v]$; that is, $\vert (N_G[u] \cup N_G[v]) \cap D\vert \geq 1$. Intuitively, a vertex-edge dominating set ensures that every edge, as well as all edges incident to either of its endpoints, is dominated by at least one vertex from the set. The \textsc{Min-VEDS} problem asks for a vertex-edge dominating set of minimum size in a given graph. This problem is known to be NP-complete even for bipartite graphs. In this paper, we strengthen this hardness result by proving that the problem remains NP-complete for two specific subclasses of bipartite graphs: star-convex and comb-convex bipartite graphs. For a graph $G$ on $n$ vertices, it is known that the \textsc{Min-VEDS} problem cannot be approximated within a factor of $(1 - \epsilon)\ln |V|$ for any $\epsilon > 0$, unless $\text{NP} \subseteq \text{DTIME}(|V|^{O(\log \log |V|)})$. We also prove that this inapproximability result holds even for star-convex and comb-convex bipartite graphs. On the positive side, we present a polynomial-time algorithm for computing a minimum vertex-edge dominating set in convex bipartite graphs. A polynomial-time algorithm for this graph class was also proposed by B{\"u}y{\"u}k{\c{c}}olak et al.~\cite{buyukccolak2025linear}, but we show that their algorithm has certain flaws by providing instances where it fails to produce an optimal solution. We address this issue by presenting a modified algorithm that correctly computes an optimal solution.
\end{abstract}


\section{Introduction}

Given a graph $G=(V,E)$, a subset $D \subseteq V$ of a graph $G = (V, E)$ is called a \emph{dominating set} of $G$ if $|N_G[v]\cap D|\geq 1$ for all $v\in V$. The \emph{domination number} is the minimum cardinality of a dominating set of $G$ and it is denoted by $\gamma(G)$. Generally, a vertex $v$ dominates itself and all of its adjacent vertices. Various variations of this classical domination problem have been explored in the literature, depending on the specific nature of the dominating power of a vertex. For further information on domination and its variants, readers are encouraged to refer to \cite{haynes1998dominationAdvanced, haynes1998fundamentals}.

In the classical domination problem, a vertex $v$ dominates itself and all of its neighbouring vertices. In contrast, another variation focuses solely on edge domination, where a vertex $v$ is considered to dominate all edges incident to any vertex in $N_G[v]$. In the literature, this variation is known as vertex-edge domination (or ve-domination). Let $e=uv\in E(G)$ be an edge, we say $e = uv \in E(G)$ is \emph{vertex-edge dominated} (or \emph{ve-dominated}) by a vertex $w$ if $w \in N_G[u] \cup N_G[v]$. A vertex-edge dominating set, or VED-set for short, of $G$ is a subset $D\subseteq V$ such that for every $uv \in E$, $|(N_G[u] \cup N_G[v]) \cap D| \geq 1$. The minimum \emph{ve-domination number} of $G$, denoted as $\gamma_{ve}(G)$, is the minimum cardinality of a ve-dominating set of a graph $G$. Furthermore, the \textsc{Min-VEDS} Problem is to find a vertex-edge dominating set of minimum cardinality for the input graph. In \cite{peters1986domination}, Peters introduced this problem in his PhD thesis. The optimization version of the problem is defined as follows:  

\medskip
\noindent\underline{\textsc{Minimum Vertex-Edge Domination} problem (Min-VEDS)}
\begin{enumerate}
	\item[] \textbf{Instance}: A graph $G=(V,E)$.
	\item[] \textbf{Solution}: $\gamma_{ve}(G)$.
\end{enumerate}


\subsection{Definitions and notations}
Let $G = (V, E)$ be a simple and undirected graph with vertex set $V(G)$ and edge set $E(G)$, when the context is clear, we use $V$ and $E$ instead of $V(G)$ and $E(G)$. The \emph{open neighbourhood} of a vertex $x\in V$ is the set of vertices $y$ that are adjacent to $x$, and it is denoted by $N_G(x)$. The \emph{closed neighborhood} of a vertex $x\in V$, denoted as $N_G[x]$, is defined by $N_G[x]=N_G(x)\cup \{x\}$.

A graph $G = (V, E)$ is said to be bipartite if $V(G)$ can be partitioned into two disjoint sets $X$ and $Y$ such that every edge joins a vertex in $X$ to a vertex in $Y$. A partition $(X, Y)$ of $V$ is called a bipartition. A bipartite graph with bipartition $(X, Y)$ of $V$ is denoted by $G = (X\cup Y, E)$. Let $n$ and $m$ denote the number of vertices and the number of edges of $G$, respectively. For a bipartite graph $G = (X\cup Y, E)$, let $|X| = n_1$, $|Y| = n_2$, and $n = n_1 + n_2$. 

A \emph{star} is a tree that contains exactly one non-pendant vertex. A \emph{comb} is a tree formed by attaching a single pendant edge to each vertex of a path. In a comb, the path is referred to as the \emph{backbone}, while the pendant vertices are called the \emph{teeth}. A bipartite graph $G = (X, Y, E)$ is called a \emph{tree-convex bipartite graph} if there exists a tree $T = (X, E_X)$ such that, for every vertex $y \in Y$, the neighbourhood of $y$ in $G$ induces a subtree of $T$. Tree-convex bipartite graphs can be recognized in linear time, and the corresponding tree $T$ can also be constructed in linear time \cite{sheng2012review}.  When $T$ is a star (resp., a comb), $G$ is called a \emph{star-convex} (resp., \emph{comb-convex}) \emph{bipartite graph}.

A bipartite graph $G=(X\cup Y,E)$ is called a \textit{bipartite chain graph} if there exist orderings of vertices of $X$ and $Y$, say $\sigma_X= (x_1,x_2, \ldots ,x_{n_1})$ and $\sigma_Y=(y_1,y_2, \ldots ,y_{n_2})$, such that $N(x_{n_1})\subseteq N(x_{n_1-1})\subseteq \ldots \subseteq N(x_2)\subseteq N(x_1)$ and $N(y_{n_2})\subseteq N(y_{n_2-1})\subseteq \ldots \subseteq N(y_2)\subseteq N(y_1)$. Such orderings of $X$ and $Y$ are called  \emph{chain orderings}, and they can be computed in linear time \cite{heggernes2007linear}.

A bipartite graph $G = (X\cup Y, E)$ is called a convex bipartite graph if the vertices of $Y$ can be ordered in such a way that for every $x \in X$, the neighbours of $x$ appear consecutively in the ordering of $Y$. In this context, $Y$ is said to exhibit the property of convexity. Such an ordering is called a convex ordering of $G$. Booth and Lueker \cite{booth1976testing} described an algorithm that determines whether a given bipartite graph $G = (X\cup Y, E)$ has the convexity property on $Y$ and, if so, generates a convex ordering of $G$. Booth and Lueker’s algorithm has a linear time complexity of $O(n + m)$. Thus, by \cite{booth1976testing}, a convex graph represented as a bipartite graph without a convex ordering can be transformed into convex form in linear time $O(n + m)$, and a convex ordering of a convex bipartite graph $G$ can be computed in linear time $O(n + m)$.

Let $G$ be a convex bipartite graph and $\sigma=\{x_1, x_2, \ldots, x_{n_1}, y_1, y_2, \ldots, y_{n_2}\}$ be  a convex ordering of $G$ with respect to $Y=\{y_1, y_2, \ldots, y_{n_2}\}$. For a vertex $x\in X$, we define $left(x)$ to be the minimum index vertex in the neighbour of $x$ with respect to $\sigma$ and $right(x)$ to be the maximum index vertex in the neighbour of $x$ with respect to $\sigma$. Similarly, we can define $left(y)$ and $right(y)$ for a vertex $y\in Y$. For two vertices $x_i, x_j\in X$ with $i<j$, let  $N_G(x_i)=\{y_{r_i}, y_{r_i+1}, \ldots, y_{s_i} \}$ and $N_G(x_j)=\{y_{r_j}, y_{r_j+1}, \ldots, y_{s_j} \}$. We write $right(x_i)\prec right(x_j)$ (resp. $left(x_i)\prec left(x_j)$) or $right(x_i)\preceq right(x_j)$ (resp. $left(x_i)\preceq left(x_j)$) if $s_i<s_j$ (resp. $r_i<r_j$) or $s_i\leq s_j$ (resp. $r_i\leq r_j$). 

An ordering $\beta = (x_1, x_2, \ldots, x_{n_x}, y_1, y_2, \ldots, y_{n_y})$ of a convex bipartite graph $G = (X\cup Y, E)$ is called a lexicographic convex ordering or simply lex-convex if $\beta$ is a convex ordering of $Y$ and if for every pair of vertices $x_i, x_j\in X$ with $i < j$, $\beta $ satisfies one of the following two conditions:
\begin{enumerate}
	\item $left(x_i)\prec left(x_j)$
	\item $left(x_i)= left(x_j)$ and $right(x_i) \preceq right(x_j)$
\end{enumerate}

From \cite{panda2013linear}, for a given convex ordering $\sigma = (x_1, x_2, \ldots, x_{n_x}, y_1, y_2, \ldots, y_{n_y})$ of $G$, a lex-convex ordering $\beta$ of $G$ can be computed in $O(n + m)$ time using a bucket sort.

In \cite{buyukccolak2025linear}, B{\"u}y{\"u}k{\c{c}}olak gave a chain decomposition of a convex bipartite graph. Let $G=(X\cup Y, E)$ be a convex bipartite graph with lex-convex ordering $\sigma=\{x_1, x_2, \ldots, x_{n_1}, y_1, y_2, \ldots, y_{n_2}\}$ and $x_{t_1}=right(y_1)$. The chain decomposition of a connected convex bipartite graph $G$ is defined as follows:

\begin{enumerate}
	\item Let $H_1$ be the chain subgraph induced by $N(y_1)\cup N(x_{t_1})$.
	
	\item Remove $H_1$ form $G$. Let $J_1$ be the set of isolated vertices in $X$ within $G[V(G)\setminus V(H_1)]$.
	
	\item Remove $J_1$ from $G$ as well.
	
	\item Repeat the steps $1$, $2$ and $3$.
	
\end{enumerate}

The decomposition ends with either a chain graph or a collection of isolated vertices. The sequence of chain graphs $H_1, H_2, \ldots$ along with the sets $J_1, J_2, \ldots$ of isolated vertices is called the chain decomposition of $G$. From now on, the chain decomposition will be denoted as $(H_1, J_1, H_2, J_2,\ldots)$. Note that we do not consider an isolated as part of a chain graph. Hence, each chain graph $H_i$ satisfies $H_i\cap X\neq \emptyset$ and $H_i\cap Y\neq \emptyset$.

\begin{lemma}[\cite{buyukccolak2025linear}]
	Let $H_1, H_2, \ldots, H_t$ with $J_1, J_2, \ldots, J_t$ be the chain decomposition of a connected convex bipartite graph $G = (X\cup Y, E)$. Then,
	\begin{enumerate}
		\item $H_i$ overlaps with $J_{i}$ and $H_{i+1}$ for $1 \leq i \leq t-1$, and
		\item $H_i$ does not overlap with $J_{i+1}$, and hence $H_{i+2}$, for $1 \leq i \leq t-2$.
	\end{enumerate}
\end{lemma}
Moreover, for $1\leq i\leq k$, each vertex in $J_i$ is adjacent to a vertex in $H_i$, and some leftmost vertices in $X$ of $H_{i+1}$ are adjacent to the rightmost vertex in $Y$ of $H_i$. Given a convex bipartite graph $G=(X\cup Y, E)$, the chain decomposition of $G$ can be obtained in $O(n)$ time \cite{buyukccolak2025linear}.

We use $[t]$ to indicate the set $\{1, 2, \ldots, t\}$ and $[t, t']$ to indicate $\{t, t+1, \ldots, t'-1, t'\}$ for $t<t'$. Let $G=(X\cup Y, E)$ be  a convex bipartite graph with lex-convex ordering, $\sigma=(x_1, x_2, \ldots, x_{n_1}, y_1, y_2, \ldots, y_{n_2})$. Also, let the sequence of chain graphs $H_1, H_2, \ldots$ along with the sets $J_1, J_2, \ldots$ of isolated vertices be the chain decomposition of $G$. For each, $H_i$ we denote $X_{H_i}=X\cap H_i$ and $Y_{H_i}=Y\cap H_i$. For a given set $X=\{x_1, x_2, \ldots, x_{n_1}\}$, we denote $X_i=\{x_i, x_{i+1}, \ldots, x_{n_1}\}$ for all $i\in [n_1]$ and $X_{[r, r']}=\{x_r, x_{r+1}, \ldots, x_{r'-1}, x_r'\}$. Similarly, for $Y=\{y_1, y_2, \ldots, y_{n_2}\}$, we denote $Y_i=\{y_i, y_{i+1}, \ldots, y_{n_2}\}$ for all $i\in [n_2]$ and $Y_{[s, s']}=\{y_s, y_{s+1}, \ldots, y_{s'-1}, y_s'\}$. 

\subsection{Existing Literature}
The \textsc{Min-VEDS} problem was first introduced by Peters in 1986 in his PhD thesis \cite{peters1986domination}. Later, Lewis extended the study by providing lower bounds for $\gamma_{ve}(G)$ for various graph classes, like $k$-regular graphs and cubic graphs \cite{lewis2007vertex}. From an algorithmic perspective, Lewis showed that the \textsc{Min-VEDS} problem is NP-complete for bipartite, chordal, planar and circle graphs. He also proposed an algorithm to solve the \textsc{Min-VEDS} problem for trees \cite{lewis2007vertex}, which was later shown to be wrong by Paul et al. \cite{paul2022results}. Paul et al. also produced efficient algorithms to solve the \textsc{Min-VEDS} problem for block graphs (a generalization of trees)\cite{paul2022results}, interval graphs and bipartite permutation graphs \cite{paul2021interval}. Very recently, B{\"u}y{\"u}k{\c{c}}olak proposed an algorithm to solve the \textsc{Min-VEDS} problem for convex bipartite graphs \cite{buyukccolak2025linear} (which generalizes the result for bipartite permutation graphs in \cite{paul2021interval}), but it seems like there is a certain flaw in the algorithm, for which it fails to produce the optimal solution for certain instances. For more algorithmic and combinatorial results, we refer to \cite{boutrig2016vertex, chitra2012global, jena2022vertex, krishnakumari2014bounds, paul2022results, zylinski2019vertex}.
\subsection{Our Results}

The subsequent sections of this manuscript are organized in the following manner:
\begin{itemize}
	\item In Section \ref{sec: convex}, first we produce an example for which the algorithm of B{\"u}y{\"u}k{\c{c}}olak in \cite{buyukccolak2025linear} is not working. Then we propose our modified algorithm, which solves the \textsc{Min-VEDS} problem for convex bipartite graphs.
	
	\item In Section \ref{sec:hardness}, we show the NP-hardness of the \textsc{Min-VEDS} problem for two very well-known subclasses of bipartite graphs: comb convex bipartite graphs and star convex bipartite graphs. We also present some additional inapproximability results for these two graph classes.
	\item Section \ref{sec:conclusion} concludes the manuscript.
	
\end{itemize}

\section{Convex bipartite graphs}\label{sec: convex}
In this section, we first present a counterexample in which the algorithm proposed by B{\"u}y{\"u}k{\c{c}}olak in \cite{buyukccolak2025linear} to compute the vertex-edge domination number for convex bipartite graphs, fails to produce an optimal solution, thereby revealing an error in their algorithm. We then introduce our modified algorithm for computing the vertex–edge domination number in convex bipartite graphs.
\subsection{Counterexample}
Before constructing the counterexample, we first introduce the algorithm given in \cite{buyukccolak2025linear}. Let $G=(X\cup Y, E)$ be a connected convex bipartite graph with lex-convex ordering $\sigma=(x_1, x_2, \ldots, x_{n_1}, y_1, y_2, \ldots, y_{n_2})$, and let $(H_1,J_1,H_2,J_2,\ldots)$ be the chain decomposition of $G$. Furthermore, assume $x_{t_i}$ is the rightmost neighbour of the leftmost vertex in $Y_{H_i}$ of $H_i$ with respect to lex-ordering $\sigma$. Then the following theorem is given in \cite{buyukccolak2025linear} to calculate the minimum VED-set of $G$.

\begin{theorem}[\cite{buyukccolak2025linear}]\label{Theo_convex_1}
	The set of vertices $x_{t_1}, x_{t_2}, \ldots$ forms a minimum VED-set of $G$.
\end{theorem}

\begin{figure}[htbp!]
	\centering
	\includegraphics[scale=0.85]{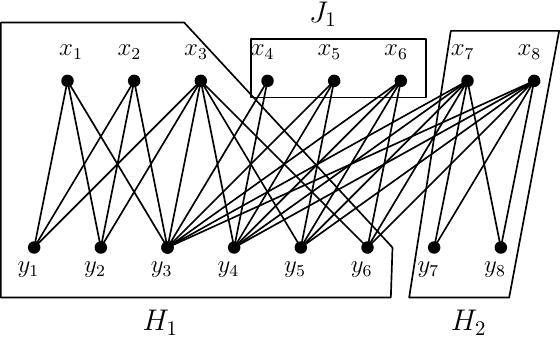}
	\caption{A convex bipartite graph with a lex-convex ordering and chain decomposition}
	\label{fig:counter_example}
\end{figure}

Now, let us consider the convex bipartite graph $G=(X\cup Y, E)$ given in Figure \ref{fig:counter_example}. This figure also presents the chain decomposition of $G$. Let $D=\{x_3, x_8\}$; according to Theorem \ref{Theo_convex_1}, the set forms a minimum VED-set. But if we take $D'=\{y_3\}$, clearly it will form a minimum VED-set of $G$. In the next section, we propose an algorithm to compute the minimum VED-set of a convex bipartite graph.

\subsection{Our Algorithm}

In this subsection, we propose a polynomial-time algorithm to compute the vertex edge domination number for convex bipartite graphs. Throughout this section, we assume that $G=(X\cup Y, E)$ is a connected convex bipartite graph, with $Y$ having the convex property.

\begin{lemma}\label{VE_convex_lemma_1}
	Let $G=(X\cup Y, E)$ be a connected convex bipartite graph with lex-convex ordering, $\sigma=(x_1, x_2, \ldots, x_{n_1}, y_1, y_2, \ldots, y_{n_2})$, and let $(H_1,J_1,H_2,J_2,\ldots)$ be the chain decomposition of $G$. Also, assume $y_{r'}= right(x_1)$, $x_r= right(y_1)$,  $y_s= right(x_r)$, $x_l=left(y_{s+1})$ and $x_p=right(y_{s+1})$. Then one of the following holds:
	\begin{enumerate}[label=\textup{(\alph*)}]
		\item There exists a minimum VED-set $D$ of $G$ such that $x_r\in D$.
		
		\item There exists a minimum VED-set $D$ of $G$ such that $y_{\alpha}\in D$ for $y_{\alpha}\in N_G(x_1)$, where $\alpha=\text{max}\{i~|~ 1\leq i\leq r' \text{ and } y_i$ is adjacent to every vertex of $J_1\}$ 
	\end{enumerate}	
\end{lemma}
\begin{proof}
	Since $G$ is a connected convex bipartite graph and $\sigma=(x_1, x_2, \ldots, x_{n_1}, y_1, y_2, \ldots, y_{n_2})$ is a lex-convex ordering, $x_1y_1\in E(G)$. Let $D$ be a minimum VED-set of $G$. In order to ve-dominate the edge $x_1y_1$, $D$ must contain at least one vertex from the set $N_G[x_1]\cup N_G[y_1]$. Next, we divide our proof into the following two cases:
	
	\begin{case}
		There exists  $x_{\beta}\in N_G(y_1)$ such that $x_{\beta}\in D$.
	\end{case}
	In this case, if $x_{\beta}=x_r$, then we are done. So let us assume $x_{\beta}\neq x_r$ and $D'=(D \setminus \{x_{\beta}\})\cup \{x_r\}$. We show that $D'$ is a minimum VED-set of $G$. Note that the vertex $x_r$ ve-dominates all the edges incident on the vertices of $\{x_r, y_1, y_2, \ldots, y_s\}$. Since $left(x_{\beta})=left(x_r)$ and $\sigma$ is a lex-convex ordering, $right(x_{\beta})\preceq right(x_r)$. It implies that $N_G(x_{\beta})\subseteq N_G(x_r)$ and hence all edges ve-dominated by the vertex $x_{\beta}$ are also ve-dominated by $x_r$. So, $D'$ is a minimum VED-set of $G$ such that $x_r\in D'$.

	\begin{case}
		There exists $y_{\alpha}\in N_G(x_1)$ such that $y_{\alpha}\in D$.
	\end{case}
	In this case, if $x_p\in D$, then $D'=(D \setminus \{y_{\alpha}\})\cup \{x_r\}$ will be our required minimum ve-dominating set. So, assume $D$ is a minimum ve-dominating set such that $x_p\notin D$. 
	
	\medskip
	\noindent \textbf{Case 2.1}. \textit{There exists a vertex $x_{a}\in J_1$ such that $x_{a}$ is not adjacent to $y_{\alpha}$}.
	
	\medskip
	Since $G$ is connected, $x_{a}$ must be adjacent to some $y_{b}\in N_G(x_r)$, but $y_{\alpha}$ is not ve-dominating the edge $x_{a}y_b$. So in order to ve-dominate the edge $x_{a}y_b$, $D$ must contain at least one vertex from the set $N_G[x_a]\cup N_G[y_b]$. Note that $N_G[x_a]\cup N_G[y_b]\subseteq X_{H_1}\cup J_1\cup Y_{H_1}\cup X_{H_2}$. Let $v\in N_G[x_a]\cup N_G[y_b]\cap D$ be the vertex in $D$ which ve-dominates the edge $x_{a}y_b$. Assume $D'=(D \setminus \{y_{\alpha}, v\})\cup \{x_p, x_r\}$. Clearly, any edge ve-dominated by $\{v, y_{\alpha}\}$ is also ve-dominated by  $\{x_p, x_r\}$. Hence, $D'$ is a minimum ve-dominating set of $G$ such that $x_r\in D'$. So, if there exists a vertex $x_{a}\in J_1$ such that $x_{a}$ is not adjacent to $y_{\alpha}$, we always have a minimum ve-dominating set of $G$ containing the vertex $x_r$. 
	
	\medskip
	\noindent \textbf{Case 2.2}. \textit{$y_{\alpha}$ is adjacent to every vertex in $J_1$}.
	
	\medskip
	In this case, if $\alpha$ is not maximum index in $[1, r']$, let $\alpha'\in [1, r']$ be the maximum index such that $y_\alpha$ is adjacent to every vertex of $J_1$. Clearly, $N_G(y_\alpha)\subseteq N_G(y_{\alpha'})$. So, $y_\alpha$ can be replaced by $y_{\alpha'}$ in $D$, which implies that $D'=(D\setminus \{y_{\alpha}\})\cup \{y_{\alpha'}\}$ is a minimum ve-dominating set of $G$ which contains $y_{\alpha'}$. 
\end{proof}

\begin{lemma}\label{VE_convex_lemma_2}
	Let $G=(X\cup Y, E)$ be a connected convex bipartite graph with lex-convex ordering, $\sigma=(x_1, x_2, \ldots, x_{n_1}, y_1, y_2, \ldots, y_{n_2})$, and let $(H_1,J_1,H_2,J_2,\ldots)$ be the chain decomposition of $G$. Also, assume $x_r= right(y_1)$, $y_s= right(x_r)$, $x_l=left(y_{s+1})$, $x_p=right(y_{s+1})$ and for some $y_{\alpha}\in N(x_1)$ assume $x_{k_{\alpha}}= right(y_{\alpha})$, $y_{l_{\alpha}}= left (x_{k_{\alpha}+1})$. Then one of the following holds.
	\begin{enumerate}[label=(\alph*)]
		\item $\gamma_{ve}(G)=\gamma_{ve}(G')+1$, where $G'=G[X_l\cup Y_{s+1}]$.
		
		\item $\gamma_{ve}(G)=\gamma_{ve}(\widetilde{G})+1$, where $\widetilde{G}=G[X_{k_{\alpha}+1}\cup Y_{l_{\alpha}}]$ for some $y_{\alpha}\in N_G(x_1)$.
	\end{enumerate}	
\end{lemma}
\begin{proof}
	Let $D'$ be a minimum VED-set of $G'$. Clearly, $D=D'\cup \{x_r\}$ is a VED-set of $G$, which implies that $\gamma_{ve}(G)\leq \gamma_{ve}(G')+1$. Similarly, if $\widetilde{D}$ is a minimum VED-set of $\widetilde{G}$ then $D=\widetilde{D}\cup \{y_{\alpha}\}$ is a VED-set of $G$ and we have $\gamma_{ve}(G)\leq \gamma_{ve}(\widetilde{G})+1$. Let $D$ be a minimum VED-set of $G$. By Lemma \ref{VE_convex_lemma_1}, we may assume that either $x_r\in D$ or there exists $y_{\alpha}\in N_G(x_1)$ such that $y_{\alpha}\in D$ and $y_{\alpha}$ is adjacent to every vertex of $J_1$. Now we divide our proof into the following two cases:
	
	\begin{case}\label{CBG_CA_1}
		Let $D$ be a minimum VED-set of $G$ such that $x_r\in D$
	\end{case}
	In this case, first we show that there is a minimum VED-set of $G$, say $\widehat{D}$, such that $x_r\in \widehat{D}$ and $X_{[1, r-1]}\cap \widehat{D}=\emptyset, Y_{[1, s]}\cap \widehat{D}=\emptyset$. If $x_p\in D$, then $D$ cannot contain any vertices from $X_{[1, r-1]}$ and $Y_{[1, s]}$. So, in this case $\widehat{D}:=D$. The other case is when $x_p\notin D$. If $y_a\in Y_{[1, s]}\cap D$, then we define $\widehat{D}:=(D\setminus {y_a})\cup \{x_p\}$. Note that no vertices from $X_{[1, r-1]}$ can be in $D$ as $x_r\in D$. So, in this case we can always assume $D$ is a minimum VED-set of $G$ such that $x_r\in D$ and $D$ does not contain any vertices from $X_{[1, r-1]}$ and $Y_{[1, s]}$.

	Since $x_r\in D$, all the edges incident with $N_G[x_r]$ are ve-dominated by $x_r$. It follows that all the edges incident on vertices $Y_{[1, s]}$ are ve-dominated by the vertex $x_r$. Now for all $i<r$, since $left(x_i)=left(x_r)$ and $\sigma=(x_1, x_2, \ldots, x_{n_1}, y_1, y_2, \ldots, y_{n_2})$ is a lex-convex ordering, $right(x_i)\preceq right(x_r)$. It implies that $N_G(x_i)\subseteq N_G(x_r)$ and hence all the edges incident to the vertices of $X_{[1, l-1]}$ are ve-dominated by $x_r$. Note that $x_r$ does not ve-dominate any edge both of whose end points lie in $J_i$ (for $i\geq 2$). Hence, every edge in $G[X_l\cup Y_{s+1}]$ is ve-dominated by $D\setminus \{x_r\}$. Let $D' = D \setminus \{x_r\}$. Thus, it is clear that $D'$ is a VED-set of $G'$. Therefore, we have $\gamma_{ve}(G') \leq |D| - 1 = \gamma_{ve}(G) - 1$. Combining this inequality with $\gamma_{ve}(G)\leq \gamma_{ve}(G')+1$, we have $\gamma_{ve}(G) = \gamma_{ve}(G') + 1$.
	

	\begin{case}
		Let $D$ be a minimum VED-set of $G$ such that $y_{\alpha}\in D$ and $y_{\alpha}$ is adjacent to every vertex of $J_1$ for some $y_{\alpha}\in N_G(x_1)$
	\end{case}
	
	First, we show that $D$ does not contain $x_p$. Otherwise, we have $D'=(D\setminus {y_{\alpha}})\cup \{x_r\}$ (note that $D'$ contains both $x_r$ and $x_p$) and we have the previous case. Similarly, we can show that if $D$ contains any vertex $v$ from $X_{[1, k_{\alpha}]}$ or $Y_{[1, \alpha -1]}$, then $D'=(D\setminus\{y_{\alpha},v\}) \cup \{x_r,x_p\}$ is a minimum VED-set which contains $x_r$, which can be handled like Case \ref{CBG_CA_1}. So, in this case, it can be assumed that $D$ is a minimum VED-set of $G$ such that $y_{\alpha}\in D$ for some $y_{\alpha}\in N_G(x_1)$ and  $y_{\alpha}$ is adjacent to every vertex of $J_1$; and $D$ does not contain any vertex from $X_{[1, k_{\alpha}]}$ and $Y_{[1, \alpha-1]}$ (for a better understanding refer to Figure \ref{fig:Y_alpha}). 
	
	\begin{figure}[htbp!]
		\centering
		\includegraphics[width=0.5\linewidth]{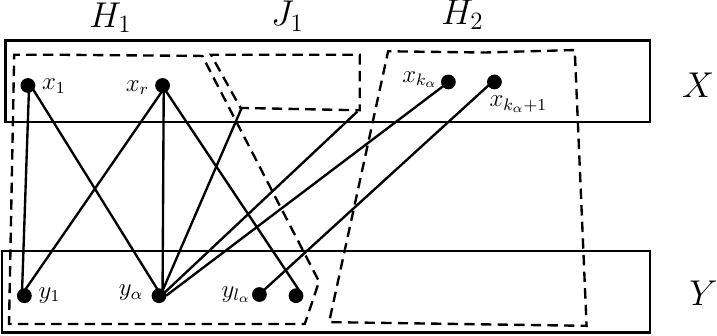}
		\caption{When $D$ contains $y_{\alpha}$}
		\label{fig:Y_alpha}
	\end{figure}
	
	Since $y_{\alpha}\in D$, all the edges incident with $N_G[y_{\alpha}]$ are ve-dominated by $y_{\alpha}$. As $y_{\alpha}$ is adjacent to every vertex of $J_1$, it follows that all the edges incident to vertices of $X_{[1, k_{\alpha}]}$ are ve-dominated by the vertex $y_{\alpha}$. Now note that $l_{\alpha}>\alpha$, or $x_{k_\alpha+1}$ is adjacent to $y_{\alpha}$, which is a contradiction. Also note that $N(y_j)\subseteq N(y_{\alpha})$ for every $j<l_{\alpha}$, and all the edges incident to the vertices of $Y_{[1, l_{\alpha}-1]}$ are ve-dominated by $y_{\alpha}$. Let $\widetilde{D}=D \setminus \{y_{\alpha}\}$. Thus, it is clear that $\widetilde{D}$ is a VED-set of $\widetilde{G}=G[X_{k_{\alpha}+1}\cup Y_{l_{\alpha}}]$.  Therefore, we have $\gamma_{ve}(\widetilde{G}) \leq |D| - 1 = \gamma_{ve}(G) - 1$. Combining this equation with $\gamma_{ve}(G)\leq \gamma_{ve}(\widetilde{G})+1$, we have $\gamma_{ve}(G) = \gamma_{ve}(\widetilde{G}) + 1$.
\end{proof}

\begin{algorithm}[htbp!]
	\textbf{Input:} A convex bipartite graph $G=(X\cup Y,E)$ with $Y$ having the convexity property. \\
	\textbf{Output:} $\gamma_{ve}(G)$.\\

	Compute the lexiographic convex ordering $\beta = (x_1, x_2, \ldots, x_{n_x}, y_1, y_2, \ldots, y_{n_y})$ and chain decomposition $(H_1,J_1,\ldots)$ of $G$;\\
	
	\uIf{$\exists~x\in X$ (or $y\in Y$) such that $N(x)=Y$ (or $N(y)=X$)}{
		\Return 1;
	}\Else{
		\uIf{$\exists ~y\in N(x_1)$ such that $y$ is adjacent to every vertex of $J_1$}{
			$\alpha=\text{max}\{i~|~ 1\leq i\leq r' \text{ and } y_i$ is adjacent to every vertex of $J_1\}$;\\
			$G'=G[X_l\cup Y_{s+1}]$;\\
			$\widetilde{G}=G[X_{k_{\alpha}+1}\cup Y_{l_{\alpha}}]$;\\
			\Return $1+$ min$\{$ALG$\_$VE$\_$DOM($G'$), ALG$\_$VE$\_$DOM($\widetilde{G}$)$\}$;
		}\Else{
			$G'=G[X_l\cup Y_{s+1}]$;\\
			\Return $1+$ min$\{$ALG$\_$VE$\_$DOM($G'$)$\}$;
		}
	}

	\caption{Algorithm to compute $\gamma_{ve}(G)$ for a connected convex bipartite graph $G$ (\textsc{ALG\_VE\_DOM($G$)})}\label{Alg:ve}
\end{algorithm}

The notations used in Algorithm \ref{Alg:ve} are directly adapted from Lemma \ref{VE_convex_lemma_2}. So, the algorithm ALG$\_$VE$\_$DOM($G$) produces $\gamma_{ve}(G)$ for a convex bipartite graph $G$. The correctness of the algorithm follows directly from Lemma \ref{VE_convex_lemma_2}. Note that the time complexity of the Algorithm \ref{Alg:ve} is $O(n^2)$.

\section{NP-completeness}\label{sec:hardness}
In this section, we show that the decision version of \textsc{Min-VEDS} is NP-complete for the star-convex bipartite graph and comb-convex bipartite graphs.

\subsection{Star-convex bipartite graphs}
Clearly, the minimum vertex-edge domination problem is in NP. We prove the NP-hardness of this problem by showing a polynomial-time reduction from \textsc{Set Cover Decision Problem}. Let $S$ be a non-empty set and $F$ a family of subsets of $S$. In the set system $(S, F)$, a subset $C \subseteq F$ is called a cover of $S$ if every element of $S$ is contained in at least one element of $C$. The \textsc{Minimum Set Cover Problem} involves finding a cover of $S$ with the smallest cardinality for a given set system $(S, F)$. For a positive integer $k$ and a set system $(S, F)$, the \textsc{Set Cover Decision Problem} asks whether $S$ has a cover of size at most $k$ which is defined as follows:

\begin{center}
	
	\fbox{%
		\parbox{0.50\linewidth}{%
			\noindent\textsc{Set Cover Decision Problem (SCDP)}
			
			\noindent\emph{Instance:} A set system $(S, F)$, integer $k$.
			
			\noindent\emph{Question:} Does $S$ has a cover of size at most $k$?%
		}%
	}
	
\end{center}

The \textsc{SCDP} is known to be NP-complete \cite{Karp1972}. Let $(S, F)$ be an instance of \textsc{SCDP}, where $S = \{s_1, s_2, \dots, s_p\}$ and $F = \{C_1, C_2, \dots, C_q\}$, with $q \leq p$. We now construct, in polynomial time, the graph $G = (X\cup Y, E)$ as follows.

\begin{itemize}
	\item For each element $s_i$ in the set $S$, consider a vertex, $a_i$ and for each set $C_j$ in the collection $F$, add a vertex $b_j$ in $G$. Let $A=\{a_1, a_2, \ldots, a_p\}$ and $B=\{b_1, b_2, \ldots, b_q\}$. If an element $s_i$ belongs to a set $ C_j$, then add an edge between vertices $a_i$ and $b_j$.
	
	\item  For every vertex $a_i\in A$, consider a vertex $z_i$ and add the edge between $a_i$ and $z_i$ for all $i\in \{1, 2, \ldots, p\}$. Let $Z=\{z_1, z_2, \ldots, z_p\}$. 
	
	\item Consider three more vertices $u, v$ and $u'$ in $G$. Also, we add two edges $uv$ and $vu'$.

	\item For all $i\in \{1, 2, \ldots, q\}$, we add an edge between $b_i$ and $u$, that is, we make a complete bipartite graph with vertex set $B\cup \{u\}$.
\end{itemize}

\begin{figure}[htbp!]
	\centering
	\includegraphics[scale=0.85]{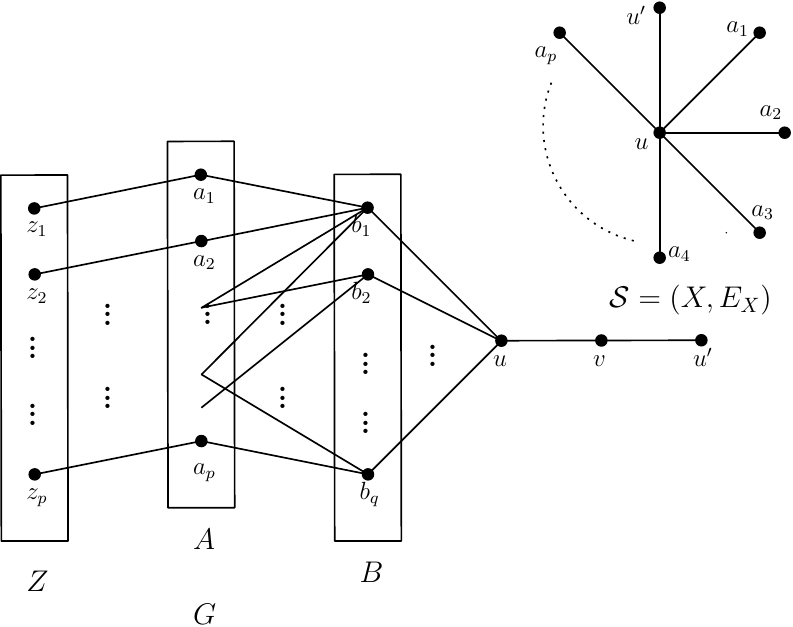}
	\caption{Construction of star-convex bipartite graph $G$ and its corresponding star $\mathcal{S}$ from a set system $(S, F)$}
	\label{fig:ve_domination_star_convex_np}
\end{figure}

Let the resultant graph be $G=(X\cup Y, E)$, where $X=A\cup \{u, u '\}$ and $Y=B\cup Z\cup \{v\}$. Clearly, $G$ is a bipartite graph. Now define a star $\mathcal{S}$ on the set $X$ with $u$ as the centre vertex. Clearly, if we take any vertex $y\in Y$, the $N(y)$ induces a subtree of $\mathcal{S}$. Hence, the constructed graph $G$ is a star-convex bipartite graph. Figure \ref{fig:ve_domination_star_convex_np} shows the construction of the star-convex bipartite graph $G$ from the set system $(S, F)$. 

Next, we show that $S$ has a cover of size $t$ if and only if $G$ has a VED-set of size at most $t+1$. For this, we have the following lemma.

\begin{lemma}\label{Lemma_VED_Star_convex}
	A set system $(S, F)$ has a cover of size $t$ if and only if $G$ has a VED-set of size at most $t+1$.
\end{lemma}

\begin{proof}
	Let $C$ be a cover of $S$ of size $t$ and $B'=\{b_j\mid C_j\in C\}$. Consider the set $D=B'\cup \{u\}$. Clearly, $|D|=t+1$. Now we prove that $D$ is a VED-set of $G$. Note that the edges $uv$ and $vu'$ are ve-dominated by the vertices $u$. Since $C$ is a cover of $S$, for every vertex $a_i\in A$, there exists $b_j\in B'$ such that $a_ib_j\in E$. So, for every edge of the form $z_ia_i$ is ve-dominated by $B'$. Furthermore, the vertex $u$ ve-dominates every edge of the form $a_ib_j$ and $ub_j$. Hence, $D$ is a VED-set of $G$ of size $t+1$.

	Conversely, let $D$ be a ve-dominating set of $G$ of size at most $t+1$. Note that to ve-dominate the edge $vu'$, $D$ must contain at least one vertex from $\{u, v, u'\}$. Clearly, when $D$ contains all $3$ vertices from $\{u, v, u'\}$, we can delete $\{v, u'\}$ from $D$ and obtain a VED-set of the size at most $t+1$. So, without loss of generality, assume $D$ contains $u$ only to ve-dominate the edge $vu'$. Let for some $i\in \{1, 2, \ldots, p\}$, $D$ contains $z_i$. In this case, if $D$ contains any vertex from the $N(a_i)$, say $b_r$, then we can remove $z_i$ from $D$ and obtain a VED-set of size at most $t+1$. On the other hand, if there are no vertices from $N(a_i)$ in $D$, we replace $z_i$ with one of the neighbours of $N(a_i)$. Similarly, if any vertex from $A$ is present in $D$, we can either remove it or replace it with its neighbour from $B$. Therefore, without loss of generality, let $D$ contains only vertices from $B$ and $\{u\}$. Let $D'=D\setminus \{u\}$. Note that $|D'|\leq t$ and $D'$ contains vertices from $B$ only. Now define $C'=\{C_i\mid b_i\in D'\}$. Let $s\in S$ and its corresponding vertex $a\in A$. We show that $a\in C_j$ for some $C_j\in C'$. Let us consider the edge $ab$, for some $b\in B$, where $b$ is the vertex corresponding to the set $C_b$ in $F$. If $C_b\in C'$, then we are done. Consider the other case when, for all $b\in N(a)$, the corresponding set $C_b\notin C'$. However, this scenario is impossible because $D$ is a VED-set that includes vertices from both $B$ and the set containing vertex $u$. Thus, $C'$ serves as a cover for $S$ with a size not exceeding $t$. Hence, $S$ has a cover of size $t$.
\end{proof}

From the above lemma, we have the following theorem.

\begin{theorem}\label{Theorem_VED_Star_convex}
	The decision version of \textsc{Min-VEDS} is NP-complete for star-convex bipartite graphs.
\end{theorem}

From Theorem \ref{Theorem_VED_Star_convex}, we have that the decision version of \textsc{Min-VEDS} is NP-complete for star-convex bipartite graphs. Next, we prove an approximation hardness result for the \textsc{Min-VEDS} in star-convex bipartite graphs. It is known that for a graph $G$, the  \textsc{Min-VEDS} cannot be approximated within a factor $(1 - \epsilon) \ln n$ for any $\epsilon>0$, unless NP $\subseteq$ DTIME($n^{O(\log \log n)}$)\cite{lewis2007vertex}. We strengthen this result for bipartite graphs by proving a similar result for star-convex bipartite graphs. For this, we need the following results by \cite{feige1995threshold} on \textsc{Minimum Set Cover Problem} and the reduction given in Theorem \ref{Theorem_VED_Star_convex}.

\begin{theorem}[\cite{feige1995threshold}] \label{TH_Min_Set_cover}
	The \textsc{Minimum Set Cover Problem} for an input instance $(S, F)$ does not admit a $(1 - \epsilon) \ln |S|$-approx algorithm for any $\epsilon > 0$ unless NP $\subseteq$ DTIME$(|S|^{O(\log \log |S|)})$. Furthermore, this inapproximability result holds for the case when the size of the input collection $F$ is no more than the size of the set $S$.
\end{theorem}

Note that the reduction given in Theorem \ref{Theorem_VED_Star_convex} is an approximation preserving reduction from \textsc{Minimum Set Cover Problem} to \textsc{Min-VEDS}. Now we have the following results.

\begin{theorem}
	The \textsc{Min-VEDS} for a star-convex bipartite graph $G$ with $|V|$ vertices does not admit a $(1 - \epsilon) \ln |V|$-approx algorithm for any $\epsilon > 0$ unless NP $\subseteq$ DTIME$(|V|^{O(\log \log |V|)})$.
\end{theorem}
\begin{proof}
	By the Lemma \ref{Lemma_VED_Star_convex}, if $D^*$ is a minimum VED set of $G$ and $C^*$ is a minimum cover of $S$ for the set system $(S, F)$, then $|D^*| = |C^*| + 1$. Now, we assume that the \textsc{Min-VEDS} for a star-convex bipartite graph can be approximated with ratio $\rho$, where $\rho = (1 - \epsilon) \ln n$ for some fixed $\epsilon > 0$, by using an approximation algorithm $\mathcal{A}$. Let $k$ be a fixed positive integer. Consider the following algorithm.
	
	\begin{algorithm}[htbp!]
		\textbf{Input:} A set $S$ and a collection $F$ of subsets of $S$. \\
		\textbf{Output:} A cover of $S$.\\

		\uIf{there exists a cover $C$ of $S$ of size $\leq k$}{
			
			$\mathcal{C}=C$;\\
			
		}\Else{ 
			Construct the corresponding graph $G$ (as explained in Theorem \ref{Theorem_VED_Star_convex});\\
			
			Compute a VED-set $D$ of $G$ using the approximation algorithm $\mathcal{A}$;\\
			
			Construct a cover $C$ of $S$ from VED-set $D$ as given in the proof of the Lemma \ref{Lemma_VED_Star_convex};\\
			
			$\mathcal{C}=C$;\\
		}	
		\Return $\mathcal{C}$;
		
		\caption{\textsc{Approx\_Set\_Cover(S, F)}}\label{Alg:SC}
	\end{algorithm}
	
	Note that the algorithm \textsc{Approx\_Set\_Cover(S, F)} is a polynomial-time algorithm, as all the steps of Algorithm \ref{Alg:SC} can be performed in polynomial time. If the cardinality of a minimum cover of $S$ is at most $k$, then it can be computed in polynomial time. Next, we analyze the case where the cardinality of the minimum cover of $S$ is greater than $k$. 
	
	Let $C^*$ denote a minimum cover of $S$ with $|C^*| > k$ and $D^*$ denote a minimum VED-set of $G$. If $\mathcal{C}$ is a cover of $S$ computed by the algorithm \textsc{Approx\_Set\_Cover(S, F)}, then 
	
	$$ |\mathcal{C}| \leq |D| \leq \rho |D^*|\leq \rho(1+ |C^*|)\leq \rho \left(1 + \frac{1}{|C^*|}\right) |C^*| < \rho \left(1 + \frac{1}{k}\right) |C^*|. $$
	
	Therefore, the algorithm \textsc{Approx\_Set\_Cover(S, F)} approximates  \textsc{Minimum Set Cover Problem} within the ratio $\rho(1 + \frac{1}{k})$. By our assumption $\rho = (1 - \epsilon) \ln |V|$ for some fixed $\epsilon > 0$. Take the integer $k> 0$ large enough so that $\frac{1}{k}< \epsilon <1$ and $\epsilon>\frac{\sqrt{3}}{2}$. Therefore, we have 
	$$|\mathcal{C}|<  (1 - \epsilon) \ln |V| \left(1+ \epsilon \right) |C^*|\leq (4-4\epsilon^2) \ln p |C^*| \quad (\text{ since } |V|\leq p^4).$$
	As $\epsilon>\frac{\sqrt{3}}{2}$, we have $\epsilon'=4\epsilon^2-3$ is a nonzero quantity and also less than $1$. Hence, $|\mathcal{C}| < (1 - \epsilon') \ln p |C^*|$, and therefore, the algorithm \textsc{Approx\_Set\_Cover(S, F)} approximates the set cover within ratio $(1 - \epsilon') \ln p$ for some $\epsilon' > 0$. By Theorem \ref{TH_Min_Set_cover}, if the \textsc{Minimum Set Cover Problem} can be approximated within $(1 - \epsilon') \ln p$, then NP $\subseteq$ DTIME$(p^{O(\log \log p)})$. This is a contradiction. Therefore, \textsc{Min-VDES} cannot be approximated within $(1 - \epsilon) \ln |V|$ for any $\epsilon > 0$ unless NP $\subseteq$ DTIME$(|V|^{O(\log \log |V|)})$.
\end{proof}

\subsection{Comb-convex bipartite graphs}

Clearly, the minimum vertex-edge domination problem is in NP. We prove the NP-hardness of this problem by showing a polynomial-time reduction from \textsc{Set Cover Decision Problem}. 

The \textsc{SCDP} is known to be NP-complete \cite{Karp1972}. Let $(S, F)$ be an instance of \textsc{SCDP}, where $S = \{s_1, s_2, \dots, s_p\}$ and $F = \{C_1, C_2, \dots, C_q\}$, with $q \leq p$. We now construct, in polynomial time, the graph $G = (X\cup Y, E)$ as follows:

\begin{itemize}
	\item For each element $s_i$ in the set $S$, consider a vertex $a_i$ and for each set $C_j$ in the collection $F$, add a vertex $b_j$ in $G$. Let $A=\{a_1, a_2, \ldots, a_p\}$ and $B=\{b_1, b_2, \ldots, b_q\}$. If an element $s_i$ belongs to set $C_j$, then add an edge between vertices $a_i$ and $b_j$.
	
	\item  For every vertex $a_i\in A$, consider a vertex $z_i$ and add the edge between $a_i$ and $z_i$ for all $i\in \{1, 2, \ldots, p\}$. Let $Z=\{z_1, z_2, \ldots, z_p\}$.
	
	\item For all $i\in \{1, 2, \ldots, p, p+1\}$, consider a vertex $r_i$ and let $R=\{r_1, r_2, \ldots, r_p, r_{p+1}\}$.
	
	\item For all $i\in \{1, 2, \ldots, q\}$ and $j\in \{1, 2, \ldots, p, p+1\}$, we add an edge between $b_i$ and $r_j$, that is, we make a complete bipartite graph with vertex set $B\cup R$.

	\item Finally, Consider two more vertices $w$ and $r_{p+1}'$ in $G$. Also, we add two edges $r_{p+1}w$ and $wr_{p+1}'$.
\end{itemize}

\begin{figure}[htbp!]
	\centering
	\includegraphics[scale=0.80]{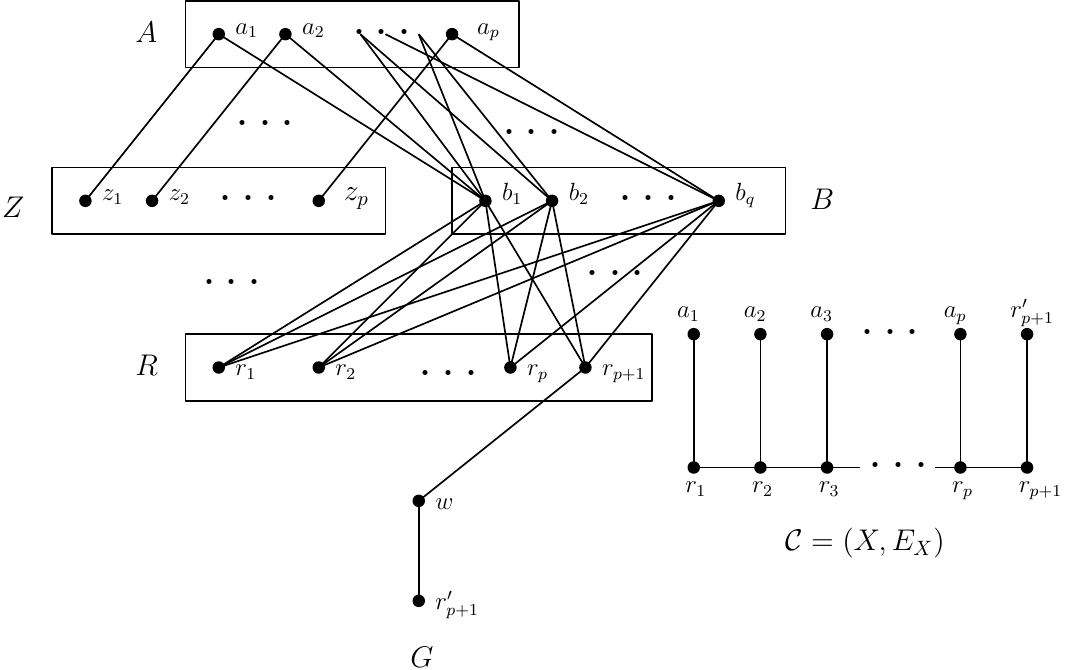}
	\caption{Construction of comb-convex bipartite graph $G$ and its corresponding comb $\mathcal{C}$ from a set system $(S, F)$}
	\label{fig:ve_domination_comb_convex_np}
\end{figure}

Let the resultant graph be $G=(X\cup Y, E)$, where $X=A\cup R\cup \{r_{p+1}'\}$ and $Y=B\cup Z\cup \{w\}$. Clearly, $G$ is a bipartite graph. Now define a comb $\mathcal{C}$ on the set $X$ with $\{r_1, r_2, \ldots, r_p, r_{p+1}\}$ as its backbone and $\{a_1, a_2, \ldots, a_p, r_{p+1}''\}$ as its teeth. Clearly, if we take any vertex $y\in Y$, the $N(y)$ induces a subtree of $\mathcal{C}$. Hence, the constructed graph $G$ is a comb-convex bipartite graph. Figure \ref{fig:ve_domination_comb_convex_np} shows the construction of the comb-convex bipartite graph $G$ from the set system $(S, F)$. 

Next, we show that $S$ has a cover of size $t$ if and only if $G$ has a VED-set of size at most $t+1$. For this, we have the following lemma.

\begin{lemma}\label{Lemma_VED_Comb_convex}
	A set system $(S, F)$ has a cover of size $t$ if and only if $G$ has a VED-set of size at most $t+1$.
\end{lemma}

\begin{proof}
	Let $C$ be a cover of $S$ of size $t$ and $B'=\{b_j\mid C_j\in C\}$. Consider the set $D=B'\cup \{r_{p+1}\}$. Clearly, $|D|=t+1$. Now we prove that $D$ is a VED-set of $G$. Note that the vertex $r_{p+1}$ ve-dominates all the edges of the form $r_lb_j$ and edges of the form $a_ib_j$. Furthermore, the edges $r_{p+1}'w$ and $wr_{p+1}$ are ve-dominated by the vertex $r_{p+1}$. Since $C$ is a cover of $S$, for every vertex $a_i\in A$, there exists $b_j\in B'$ such that $a_ib_j\in E$. Therefore, every edge of the form $z_ia_i$ is ve-dominated by $B'$. Hence, $D$ is a VED-set of $G$ of size $t+1$.

	Conversely, let $D$ be a ve-dominating set of $G$ of size at most $t+1$. Note that to ve-dominate the edge $wr_{p+1}'$, $D$ must contain at least one vertex from $\{r_{p+1}, w, r_{p+1}'\}$. Clearly, when $D$ contains all $3$ vertices from $\{r_{p+1}, w, r_{p+1}'\}$, we can delete $\{w, r_{p+1}'\}$ from $D$ and obtain a VED-set of the size at most $t+1$. So, without loss of generality, assume $D$ contains $r_{p+1}$ only to ve-dominate the edge $wr_{p+1}'$. As $D$ contains $r_{p+1}$, if $D$ contains any vertex from $\{r_1, r_2, \ldots, r_p\}$, we can remove them from $D$ and obtain a VED-set of size at most $t+1$.

	Let for some $i\in \{1, 2, \ldots, p\}$, $D$ contains $z_i$. In this case, if $D$ contains any vertex from the $N(a_i)$, say $b_r$, then we can remove $z_i$ from $D$ and obtain a VED-set of size at most $t+1$. On the other hand, if there are no vertices from $N(a_i)$ in $D$, we replace $z_i$ with one of the neighbours of $N(a_i)$. Similarly, if any vertex from $A$ is present in $D$, we can either remove it or replace it with its neighbour from $B$. Therefore, without loss of generality, let $D$ contains only vertices from $B$ and $\{r_{p+1}\}$. Let $D'=D\setminus \{r_{p+1}\}$. Note that $|D'|\leq t$ and $D'$ contains vertices from $B$ only. Now define $C'=\{C_i\mid b_i\in D'\}$. Let $s\in S$ and its corresponding vertex $a\in A$. We show that $a\in C_j$ for some $C_j\in C'$. Let us consider the edge $ab$, for some $b\in B$, where $b$ is the vertex corresponding to the set $C_b$ in $F$. If $C_b\in C'$, then we are done. Consider the other case when, for all $b\in N(a)$, the corresponding set $C_b\notin C'$. However, this scenario is impossible because $D$ is a VED-set that includes vertices from both $B$ and the set containing vertex $r_{p+1}$. Thus, $C'$ serves as a cover for $S$ with a size not exceeding $t$. Hence, $S$ has a cover of size $t$.
\end{proof}

From the above lemma, we have the following theorem.

\begin{theorem}\label{Theorem_VED_Comb_convex}
	The decision version of \textsc{Min-VEDS} is NP-complete for comb-convex bipartite graphs.
\end{theorem}

Similarly, we prove an approximation hardness result for the \textsc{Min-VEDS} in comb-convex bipartite graphs. Note that the reduction given in Theorem \ref{Theorem_VED_Comb_convex} is an approximation preserving reduction from \textsc{Minimum Set Cover Problem} to \textsc{Min-VEDS}. Now we have the following theorem for a comb-convex bipartite graph.

\begin{theorem}
	The \textsc{Min-VEDS} for a comb-convex bipartite graph $G$ with $|V|$ vertices does not admit a $(1 - \epsilon) \ln |V|$-approx algorithm for any $\epsilon > 0$ unless NP $\subseteq$ DTIME$(|V|^{O(\log \log |V|)})$.
\end{theorem}

\section{Conclusion}\label{sec:conclusion}
In this paper, we have designed a polynomial-time algorithm for computing a minimum vertex-edge dominating set of convex bipartite graphs. Furthermore, we have shown that the decision version of \textsc{Min-VEDS} remains NP-complete for both star-convex and comb-convex bipartite graphs, which are subclasses of tree-convex bipartite graphs. Additionally, we have demonstrated that for any $\epsilon > 0$, the \textsc{Min-VEDS} problem cannot be approximated within a factor of $(1 - \epsilon) \ln |V|$, even for star-convex and comb-convex bipartite graphs, unless $\mathrm{NP} \subseteq \mathrm{DTIME}(|V|^{O(\log \log |V|)})$. It would be interesting to study the complexity status of this problem in other subclasses of bipartite graphs, such as circular convex bipartite and triad convex bipartite.

\bibliographystyle{plain}
\bibliography{VE_domination_ref}

\end{document}